\documentclass[11pt]{amsart}

\usepackage{amsmath,amsfonts,amssymb,layout}
\numberwithin{equation}{section}

\def\N{{\mathcal N}}
\def\P{{\mathbb P}}
\def\Q{{\mathbb Q}}
\def\G{{\mathbb G}}
\def\C{{\mathbb C}}
\def\Z{{\mathbb Z}}
\def\Ho{\mathrm{Hodge}}
\def\Om{{\Omega}}
\def\qed{\hfill\vbox{\hrule\hbox{\vrule\kern3pt\vbox{\kern6pt}\kern3pt\vrule}\hrule}\bigskip}

\normalbaselines
\newtheorem{prop}{\bf Proposition}[section]
\newtheorem{tma}[prop]{Theorem}
\newtheorem{lma}[prop]{Lemma}
\newtheorem{cor}[prop]{Corollary}

\newtheorem{rem}[prop]{Remark}
\newtheorem{ex}[prop]{Example}
\newenvironment{pf}{{\it Proof.}}{\qed}
\newenvironment{df}{\bigskip\noindent{\bf Definition}}

\begin{document}

\title{On the Picard group of low-codimension subvarieties}
\author{Enrique Arrondo}
\author{Jorge Caravantes}
\address{ Departamento de \'Algebra\\ Facultad de Ciencias Matem\'aticas\\
Universidad Complutense de Madrid\\ 28040 Madrid, Spain}
\email{arrondo@mat.ucm.es, jcaravan@mat.ucm.es}
\thanks{Both authors were partially supported by the MCYT project
BFM2003-03971}

\begin{abstract} We introduce a method to determine if
$n$-dimensional smooth subvarieties of an ambient space of
dimension at most $2n-2$ inherit the Picard group from the
ambient space (as it happens when the ambient space is a
projective space, according to results of Barth and Larsen).
As an application, we give an affirmative answer (up to some
mild natural numerical conditions) when the ambient space is
a Grassmannian of lines (thus improving results of Barth, Van
de Ven and Sommese) or a product of two projective spaces of
the same dimension.
\end{abstract}

\maketitle

A (complex) smooth subvariety of small codimension is
expected to be very special. The most characteristic example
of this principle is the case of subvarieties of the
projective space, in which a famous conjecture of Hartshorne
(see \cite{H}) states that a smooth $n$-dimensional
subvariety $X\subset\P^N$ must be a complete intersection if
$2N<3n$. The main evidence for this conjecture is a series
of results of Barth and Larsen showing that the topology of
$X$ behaves like the one of a complete intersection.
Specifically, it was first proved in \cite{B} that
$H^i(X,\Q)\cong H^i(\P^N,\Q)$ if $i\le 2n-N$, then in
\cite{BL} that $\pi_1(X)=0$ if $N\le 2n-1$ and finally in
\cite{L} that $H^i(X,\Z)\cong H^i(\P^N,\Z)$ if $i\le 2n-N$.
In particular, this implies that the Picard group of such
$X$ is generated by the class of its hyperplane section if
$N\le 2n-2$. 

More generally, one can replace the projective space with
other ambient spaces to study this kind of problems. For
instance, in \cite{BV}, Barth and Van de Ven prove that if
$X$ is a smooth codimension-two subvariety of $\G(k,N)$ (the
Grassmann variety of subspaces of dimension $k$ in $\P^N$),
then $H^i(X,\Z)\cong H^i(G(k,N),\Z)$ if $N-k\ge6$ and
$i=0,1,2$. A stronger result has been obtained by Sommese
in \cite{S}, in the general framework of homogeneous
varieties (with particularly strong results for Abelian
varieties and products of Grassmannians). When
particularized to $X\subset Y:=\prod\G(k_l,N_l)$ of dimension $n$, his
result states that $\pi_i(Y,X)=0$ if
$i\le 2n-2\sum[(k_l+1)(N_-k_l)]+\mbox{min}\{N_l\}+1$. As a consequence, if $Y=\G(k,N)$, then the Picard
group of $X$ is generated by the hyperplane section when 
$2(k+1)(N-k)-N\le 2n-2$ (for $k=0$ this is the mentioned
result for the projective space). We can also apply this to products of projective spaces and we get that, if $Y:=\P^N\times\P^N$, then the Picard group of $X$ is generated by the restrictions of the generators of Pic$(\P^N\times\P^N)$ when $3N+2\le 2n$.

In this paper, we introduce a general method to study when
the Picard group of an ambient space $Y$ of dimension $N$ is
inherited by smooth subvarieties of dimension at least
${N+2\over2}$ (hence in the range of Barth-Larsen result for
projective spaces). Our starting point is the paper by the
first author \cite{A}, where it is shown that the canonical
divisor of any  smooth fourfold in $\G(1,4)$ (and in $Q_6$
under some numerical conditions) is --numerically-- a
multiple of the hyperplane section. The method there relied
strongly on the fact that the codimension was two (in order
to use the so-called Hartshorne-Serre correspondence) and it
seemed essential to work with the canonical divisor and not
any other; a trivial application of Porteous formula
produced a mysterious numerical relation, that eventually
implied (via the Hodge index theorem) that the surface
obtained by successive hyperplane sections was (numerically)
subcanonical. 

Since Sommese's result indicates that the result of \cite{A}
should be true in higher codimension, we tried to
understand where that mysterious numerical relation
actually came from, without using such sophisticated
techniques. We have thus found out that, in the general
framework of $X\subset Y$ of dimension $n\ge{N+2\over2}$,
that relation could be easily obtained by just applying the
self-intersection formula for $X$ in $Y$, the
self-intersection formula for any smooth divisor
$D$ of $X$ (this shows why it is important that $D$ has
dimension at least its codimension in this ambient space)
and the relation among them. In this way, we obtain in
several ambient spaces (not necessarily with Picard number
one) that any divisor (with only few restrictions) of $X$
must be numerically dependent on the generators of the
Picard group of $Y$. When $X$ is also simply connected, this
method shows that the Picard group of $X$ (which is
then torsion-free) is generated (with rational
coefficients) by the generators of the Picard group of $Y$.

Our first example, to which we devote Section \ref{proj}, is 
$Y=\P^N$ and in this case we reobtain (up to divisibility) the
theorem of Barth-Larsen for the Picard groups of subvarieties
$X\subset\P^N$. In fact, as it happened in \cite{A}, the
above mentioned numerical equality becomes
$H_{|S}^2D_{|S}^2=(H_{|S}D_{|S})^2$, where $S$ is the surface
obtained intersecting $X$ with $n-2$ general hyperplanes.
The Hodge index theorem implies immediately that $D_{|S}$ is
numerically a multiple of $H_{|S}$. The Lefschetz hyperplane
section and the simple connectedness of $X$ obtained in \cite{BL}
completes the proof.

Following the idea of the projective case, in the second
section we apply this method when the ambient space is a
Grassmannian of lines (observe that our hypothesis about the relation between the dimension and codimension of $X$ are much
milder than the ones in the result of Sommese). Here we prove
that, when $X$ has positive intersection with all Schubert
cycles of codimension $n$, then its Picard group has rank
one. In this case, the key point is to write the numerical
formula as a positive linear combination of Hodge-type
expressions for the surfaces obtained by intersecting
$X$ with the different Schubert cycles of codimension $n-2$
(which, by Hodge index theorem, must be positive, and they
are zero when numerical equivalence occurs), plus some other
positive terms. 

In the third section we illustrate the fact that our method
keep working even if the Picard group of the ambient space
$Y$ has rank bigger than one. As a sample, we have chosen
the case in which $Y$ is the cartesian product of two
projective spaces of the same dimension. This particular
choice has the advantage that $Y$ has even dimension (this
is also one of the reasons why we choose $k=1$ in the case
of Grassmannians), which simplifies the calculations. Again
we prove in this case that, in the Barth-Larsen range, the
Picard group of $X\subset Y$ is generated by two elements,
assuming again some nonvanishing of the coefficients of the
class of $X$ in the Chow ring of $Y$. This is again an
improvement of the range given in \cite{S}, where the
result is however stronger, in the sense that the Picard
group of $X$ is proved to be generated by the restriction
of the Picard group of the ambient space.

As a final section, we discuss the general validity of
the method. We first briefly treat some other examples, in
order to show that the method works also for non-homogeneous
ambient spaces, and that in some cases extra assumptions are
needed. We conclude with some general remarks and
conjectures.

\section{Subvarieties of projective space and
general method}\label{proj}

In this section we will prove (up to divisibility of the
hyperplane divisor class) Barth-Larsen theorem for the
Picard group of subvarieties of projective spaces. This will
illustrate the method that we will use later for different
ambient varieties.

\subsection{The Barth-Larsen result} 

First of all, we prove that the group of divisors up to numerical equivalence
of a low codimension variety in $\P^N$ has rank 1.

\begin{prop}\label{gorda_proj}
Let $X$ be a nonsingular projective variety of dimension $n$ in 
$\P^N$ with $N\le 2n-2$ and let $H_{|X}$ denote the
restriction to $X$ of the hyperplane class of $\P^N$. Then,
for any divisor $D$ of $X$, there exists $q\in\Q$ so that
$D\equiv_{\mathrm{num}} qH_{|X}$.
\end{prop}

\begin{pf}
Regarding $\P^N$ as a linear subspace of $\P^{2n-2}$, we can
assume $N=2n-2$. Moreover, changing $D$ with $D+kH$
for a sufficient large $k\in\Z$, we can also assume that $D$
is very ample and hence smooth. Thus, we have the following
exact sequence of vector bundles on $D$:
$$0\to \N_{D/X}\to \N_{D/\P^{2n-2}}\to (\N_{X/\P^{2n-2}})_{|D}\to 0$$

Hence we get an equality for the top Chern
classes:
\begin{equation}\label{chern_proj}
c_{n-1}(\N_{D/\P^{2n-2}})=c_1(\N_{D/X})c_{n-2}((\N_{X/\P^{2n-2}})_{|D})
\end{equation}

Let $d$ and $\delta$ be the degrees of $X$ and $D$ in
$\P^{2n-2}$  respectively. In terms of the intersection of
cycles in $X$ we have $d=H_{|X}^n$ and
$\delta=DH_{|X}^{n-1}$. On the other hand, the
self-intersection formula applied to each of the inclusions
$D\subset  X\subset \P^{2n-2}$ yields
$$c_{n-1}(\N_{D/\P^N})=\delta^2=(DH_{|X}^{n-1})^2$$ 
$$c_1(\N_{D/X})=D_{|D}$$
$$c_{n-2}(\N_{X/\P^N})=dH_{|D}^{n-2}$$

Therefore (\ref{chern_proj}) becomes:
\begin{equation}\label{autoint_proj}
(DH_{|X}^{n-1})^2=dD^2H_{|X}^{n-2}=H_{|X}^n(D^2H_{|X}^{n-2})
\end{equation}
where all the intersection products are considered in $X$.

Taking $S$ to be the smooth irreducible surface obtained as 
intersection of $X$ with $n-2$ general hyperplanes, formula 
(\ref{autoint_proj}) becomes
\begin{equation}\label{Hodge_proj}
(D_{|S}H_{|S})^2 = D_{|S}^2H_{|S}^2
\end{equation}
Thus, Hodge index theorem implies that there exists $q\in
\Q$ such that $D_{|S}\equiv_{\mathrm{num}} qH_{|S}$. 
Finally, by Lefschetz's hyperplane theorem, the
restriction map $H^2(X,\Q)\to H^2(S,\Q)$ is injective, hence
$D \equiv_{\mathrm{num}} qH_{|X}$.
\end{pf}

We obtain immediately from this the following consequence of
the results of Barth and Larsen (observe that we have to use
only their result about simple connectedness, which is in
some sense weaker):

\begin{cor}\label{Barth-Larsen_proj} In the hypothesis of
Proposition \ref{gorda_proj}, Pic$X \simeq \Z$.
\end{cor}

\begin{pf}
Proposition \ref{gorda_proj} shows that Pic$X$ has rank one,
so it is enough to show that it is torsion-free. But it is an
immediate consequence of the fact that, as proved in
\cite{BL}, $X$ is simply connected.
\end{pf}

\subsection{The general method}\label{metodo_general}

The above proofs suggest the following strategy to study
when the Picard group is inherited (up to divisibility) by
subvarieties of small codimension. Our general set-up will
be a smooth ambient space $Y$ of dimension $N\ge6$ and a
smooth subvariety $X\subset Y$ of dimension $n$, and we will
assume $N\le 2n-2$. For simplicity, assume that $N$ is even,
so that after intersecting with ${2n-2-N\over2}$ general very
ample divisors (and using Lefschetz's hyperplane theorem) we
can assume $N=2n-2$. The first step now is to take a smooth
divisor $D$ and use the exact sequence
$$0\to \N_{D/X}\to \N_{D/Y}\to (\N_{X/Y})_{|D}\to 0$$
to obtain an equality
$$c_{n-1}(\N_{D/Y})-c_1(\N_{D/X})c_{n-2}((\N_{X/Y})_{|D})=0$$
This equality becomes, using the self-intersection formulas
of the three inclusions $D\subset X\subset Y$ the equality
becomes
\begin{equation}\label{chern_general}
P:=D\cdot D-D^2\cdot X_{|X}=0,
\end{equation}
where the first summand is a product in $Y$ and the second
one is a product in $X$ (in this way, when $D=K$, this is
exactly the ``mysterious'' equality (5) of \cite{A} for a
subvariety of $\G(1,4)$). The key ingredients to prove that
the divisors on $X$ come (up to divisibility) from divisors
on $Y$ will be then the following, corresponding
respectively to Proposition
\ref{gorda_proj} and Corollary \ref{Barth-Larsen_proj} in the
projective case: 

{\it Step (i):} Derive from (\ref{chern_general}) and the
Hodge index formula on the surface $S$ (the intersection
of $n-2$ general very ample divisor on $X$) that the
restriction of $D$ to the surface $S$ is numerically
equivalent to a linear combination of restrictions of
a set of generators of Pic$Y$. This would imply, by
Lefschetz's hyperplane theorem, that any divisor on $X$ is
numerically equivalent (with rational coefficients) to the
restriction of some divisor on $Y$.

{\it Step (ii):} If $X$ is simply connected, then from (i) we get
that in fact any divisor has a multiple that is a
restriction of a divisor on $Y$.

Observe that step (i) was immediate in the case of the
projective space. In fact, equality (\ref{chern_general}) is
equivalent to say (via the Hodge index theorem) that the
restriction of $D$ to $S$ is a multiple of the hyperplane
section. In general, this will be not so easy, and we will
need to decompose $P$ as a sum of positive summands, as it
happened in \cite{A} when $Y=\G(1,4)$. Following this
example, one should expect some of the summands to be
Hodge formulas for the intersection of $X$ with suitable
intersections of $X$ with cycles in $Y$ of codimension
$n-2$ (one extra problem here is to show that these surfaces
to which apply the Hodge index theorem are smooth
irreducible). This is in fact the trickiest part of
the method, and we do not see a way to decide a priori
whether such a decomposition is possible or not. 

About step (ii), as far as we know there are no general results that lead
to expect that smooth subvarieties of small codimension are simply
connected. We thus decided to concentrate in this paper on cases for which we
now a priori that the simple connectedness hold. This is one of the reasons
for having chosen Grassmannians and products of projective spaces, since
Debarre proved in \cite{D} the kind of result we needed (also the one about
the irreducibility problem mentioned in the above paragraph).

\section{Subvarieties of Grassmannians of lines}

In this section we apply the method we explained in the
previous section to low-codimension subvarieties in
Grassmannians. We will choose Grassmannians of lines
because of two reasons. The first one is that, as we
observed, it is convenient to work on an ambient space of
even dimension. The second and main reason is that we need to
use Schubert calculus, and the use of a compact notation can be very
difficult and confusing. For instance, already for Grassmannians of lines in
$\P^n$ the rank of the Chow ring in a given codimension depends on the
parity of $n$. We thus decided to use integral parts in our expressions,
hoping that this will not make the reading obscure; a reasonable choice
could be to read this section twice: one of them thinking of
$n=2k$ and another one thinking of $n=2k+1$. For Grassmannians of higher
dimensional subspaces it seems better to study in each
case the concrete Grassmannian that one could need to
deal with.

\subsection{Statement of the result and first remarks}


Since we will need it continuously, we first recall the
following:

\begin{df} A {\it Schubert cycle} $\Omega(a,b)$ (with $0\le
a<b\le n$) is the class, in the Chow ring of $\G(1,n)$, of the
($a+b-1$)-dimensional {\it Schubert variety} $\Omega(A,B)$ consisting 
of the lines of $\P^n$ meeting the linear subspace $A$ and contained in the
linear subspace $B$, where $A\subset B\subset\P^n$, $\dim A=a$ and $\dim
B=b$. We refer to \cite{KL} for the needed background of Schubert cycles and
their intersections.
\end{df} 

The main result that we will prove is the following:

\begin{tma}
\label{Barth-Larsen_grass}
Let $X\subset\G(1,n)$ be a smooth subvariety of dimension
$n'\ge n$ and let $H$ be denote the hyperplane class of $X$
(after the Pl\"ucker embedding). If all intersections of
$X\cdot H^{n'-n}$ with the Schubert cycles of codimension
$n$ are different from zero, then Pic$X\simeq\Z$.
\end{tma}

\begin{rem}{\rm 
Let $X\subset\G(1,n)$ be a smooth irreducible 
subvariety of codimension two. It will have class
$$[X]=a_1\Omega(n-3,n)+a_2\Omega(n-2,n-1).$$
It is a simple exercise in Schubert calculus to see
that
$$[X]H=a_1\Omega(n-4,n)+(a_1+a_2)\Omega(n-3,n-1).$$
Hence, if $a_1\neq0$, $X$ is in the hypothesis of Theorem
\ref{Barth-Larsen_grass} if $n\ge5$, so that it is not
necessary to impose $a_2\ne0$. This shows that in general,
when $n'>n$, it is not necessary to impose the sufficient
condition that the intersection of $X$ with all the Schubert
cycles of codimension $n'$ are different from zero.

Moreover, in these conditions, if $a_1=0$ this means that the union in $\P^n$
of the $(2n-4)$-dimensional family of lines parametrized by
$X$ has dimension at most $n-1$. This easily implies that
$X$ is the Schubert variety of lines contained in a
hyperplane of $\P^n$ and thus its Picard group is generated
by its hyperplane section (and hence even in this case the thesis of the
theorem remains true).

Also in the case $n=4$, having in mind that the intersection
of $X$ with the Schubert cycles of lines in a
general hyperplane of $\P^4$ is smooth, it is not difficult
to check that $a_1,a_2\ne 0$ unless $X$ is the Schubert
variety of lines in a general hyperplane of $\P^4$. 

Summing up, Theorem \ref{Barth-Larsen_grass} implies that
the Picard group of any smooth subvariety
$X\subset\G(1,n)$ of codimension two is generated by one
element if $n\ge 4$. This improves Barth, van de Ven and
Sommese results in \cite{BV} and \cite{S} (which require respectively
$n\ge7$ and $n\ge 6$), although unfortunately for $n=4,5$ we are not able to
prove that we can take the hyperplane section as a generator of the Picard
group. 
}\end{rem}

\begin{rem}{\rm
In general, the numerical hypothesis in Theorem
\ref{Barth-Larsen_grass} are necessary. Consider for example
the set of all the lines of a smooth quadric in
$\P^5$. This is a smooth subvariety of $\G(1,5)$ of
codimension 3 whose intersection with the Schubert cycle of
all the lines passing through a point of $\P^5$ is clearly
zero. On the other hand, identifying the smooth quadric
with the Pl\"ucker embedding of $\G(1,3)$ in $\P^5$, it
follows that $X$ is canonically isomorphic to the incidence
variety in $\P^3\times{\P^3}^*$, whose Picard group
has rank two.
}\end{rem}

We now make some reductions to prove Theorem \ref{Barth-Larsen_grass} in the
line of subsection \ref{metodo_general}.

\begin{rem}{\rm\label{reduccion_grass}
Debarre has proved in \cite[Corollaire 7.4]{D}  that an $X$
in the hypothesis of Theorem \ref{Barth-Larsen_grass} is
simply connected. Hence, as explained in subsection
\ref{metodo_general}, we can assume $n'=n$ and if $D$ is a
smooth divisor of $X$ and $S$ is the intersection of $X$ with
$n-2$ general hyperplanes, it is enough to prove 
$D_{|S}\equiv_{num}qH_{|S}$ for some $q\in\Q$.  
}\end{rem}

We will devote the rest of the section to prove Theorem
\ref{Barth-Larsen_grass} according to Remark
\ref{reduccion_grass}.

\subsection{General set-up} Now we state the notation for
the numerical data that we will use throughout the section.
We will fix a smooth $n$-dimensional subvariety
$X\subset\G(1,n)$ with a smooth divisor $D$ on it. We will
write the classes of $X$ and $D$ in the Chow ring of
$\G(1,n)$ as
\begin{equation}\label{clase_X_en_G}
[X]=a_1\Omega(1,n)+a_2\Omega(2,n-1)+\ldots+
a_{[{n\over2}]}\Omega(\big[{n\over2}\big],\big[{n+3\over2}\big])
\end{equation}
(observe that the hypothesis about $X$ in the
statement of Theorem \ref{Barth-Larsen_grass} is saying that
all $a_i$ are different from zero)
\begin{equation}\label{clase_D_en_G}
[D]=\alpha_1\Omega(0,n)+\alpha_2\Omega(1,n-1)+\ldots+
\alpha_{[{n+1\over2}]}\Omega(\big[{n-1\over2}\big],\big[{n+2\over2}\big]).
\end{equation}

A standard Schubert calculus provides then the following
intersection table of cycles of $X$, which also defines the
numbers $\lambda_i$:
\begin{equation*}
\begin{tabular}{c|cccccccc}
\ & \scriptsize$\Om(0,n)_X$ & \scriptsize$\Om(1,n-1)_X$
& \dots &
\scriptsize$\Omega([{n-1\over2}],[{n+2\over2}])_X$
& \scriptsize$D\Om(1,n)_X$ &
\scriptsize$D\Om(2,n-1)_X$ &
\dots &
\scriptsize$D\Omega([{n\over2}],[{n+3\over2}])_X$
\\
\hline
$H_X$ & $a_1$ & $a_1+a_2$ & \dots &
$a_{[{n-1\over2}]}+a_{[{n+1\over2}]}$ &
$\alpha_1+\alpha_2$ &
$\alpha_2+\alpha_3$ & \dots &
$\alpha_{[{n\over2}]}+\alpha_{[{n+2\over2}]}$ \\
$D$ & $\alpha_1$ & $\alpha_2$ & \dots &
$\alpha_{[{n+1\over2}]}$ &
$\lambda_1$ & $\lambda_2$ & \dots & $\lambda_{[{n\over2}]}$
\end{tabular}
\end{equation*}
(where the notation $\Omega(i,j)_X$ represents the restriction to $X$ of the
Schubert cycle $\Omega(i,j)$).

\begin{rem}{\rm
Observe that, if $n$ is odd, the term
$a_{[{n+1\over2}]}$ is not defined, so that we will take it
to be zero. In the same way, we define
$\alpha_{[{n+2\over2}]}$ to be zero if $n$ is even. To have
a compact expression for the matrix we can also set
$a_0$=0. As we remarked at the beginning of the section, this compact
expression, although convenient for avoiding unnecessary repetitions, can be
confusing. Hence for the convenience of the reader we write the precise form
of the above matrix of intersections depending on whether
$n$ is even or odd:

\begin{equation}\label{matrizota_grass}
M:=\left\{\begin{matrix}
\begin{pmatrix}
a_1 & a_1+a_2 & \ldots & a_{k-1}+a_k
&
\alpha_1+\alpha_2 & \alpha_2+\alpha_3 & \ldots
&\alpha_k \\
\alpha_1 & \alpha_2 & \ldots & \alpha_k &
\lambda_1 &
\lambda_2 & \ldots &\lambda_k
\end{pmatrix}& {\rm if}\ n=2k\\ \\
\begin{pmatrix}
a_1 & a_1+a_2 & \ldots & a_k
&
\alpha_1+\alpha_2 & \alpha_2+\alpha_3 & \ldots
&\alpha_k+\alpha_{k+1} \\
\alpha_1 & \alpha_2 & \ldots & \alpha_{k+1} &
\lambda_1 &
\lambda_2 & \ldots &\lambda_k
\end{pmatrix}& {\rm if}\ n=2k+1
\end{matrix}\right.
\end{equation}

}\end{rem}

If $D$ is numerically a multiple of $H_X$, it is clear that
the matrix $M$
must have rank one. Reciprocally, we claim that if $M$ has
rank one, then Theorem \ref{Barth-Larsen_grass} follows.
Indeed if there is $q\in\Q$ such that the second row of
$M$ is $q$ times the first one, this means that the product
of $$D':=D-qH_X$$ 
with any $\Omega(i,n-i)_X$ and any $D\Omega(j,n+1-j)_X$ is
zero. Since $H_X^{n-1}$ and $H_X^{n-2}$ are respectively a
linear combination of Schubert cycles of the type
$\Omega(i,n-i)_X$ and $\Omega(j,n+1-j)_X$, it follows
that $D'H_X^{n-1}=D'DH_X^{n-2}=0$. Hence, if 
$S$ is the intersection of $X$ with $n-2$ general
hyperplanes, it follows:
$$D'_{|S}H_{|S}=D'H_{X}^{n-1}=0$$
$${D'_{|S}}^2=D'(D-qH_X)H_X^{n-2}=0$$
By Hodge index theorem, we conclude that $D'_{|S}$ is
numerically equivalent to zero. By Remark
\ref{reduccion_grass}, Theorem \ref{Barth-Larsen_grass}
follows.

We will thus devote the rest of the section to
prove that $M$ has rank one (Corollary \ref{rango-uno}).

\subsection{End of the proof}

\begin{lma}\label{lema_autoint_grass}
The following equality holds.
\begin{equation}\label{autoint_grass}
P:=\alpha_1^2+\ldots+\alpha_{[{n+1\over2}]}^2
-a_1\lambda_1-\ldots-a_{[{n\over2}]}\lambda_{[{n\over2}]}=0
\end{equation}
\end{lma}

\begin{pf}
This is (\ref{chern_general}) using (\ref{clase_X_en_G}),
(\ref{clase_D_en_G}) and the intersection products
given by Schubert calculus and the matrix $M$.
\end{pf}

\begin{rem}\label{Si_irreducibles}{\rm
For all $i=1,\ldots,[{n\over2}]$, let $S_i$ be the surface 
obtained intersecting $X$ with a general $\Om(i,n-i+1)$. It
is easy to see that $S_i$ is nonsingular, and by
\cite[Th\'eor\`eme 8.1]{D} it is also irreducible.

Applying Hodge index theorem to $S_i$ and using the fact that
$H\Omega(i,n-i+1)=\Omega(i-1,n-i+1)+\Omega(i,n-i)$, we get
\begin{equation}\label{Hodge_grass}
0\ge\Ho_i:=\begin{vmatrix} H_{|S_i}^2&D_{|S_i}H_{|S_i}\\
D_{|S_i}H_{|S_i}&D_{|S_i}^2
\end{vmatrix}=\begin{vmatrix} a_{i-1}+a_i & \alpha_i+\alpha_{i+1}\\
\alpha_i & \lambda_i \end{vmatrix}+
\begin{vmatrix} a_i+a_{i+1} & \alpha_i+\alpha_{i+1}\\
\alpha_{i+1} & \lambda_i \end{vmatrix}
\end{equation}
(as usual, one should be careful with the correct expression
of $Hodge_i$ in the cases $i=1,[{n\over2}]$; in particular,
notice that, if $n$ is even and $i=[{n\over2}]$, then the
correct definition of $Hodge_i$ consists only of the first
summand of (\ref{Hodge_grass})).
}\end{rem}

\begin{lma}\label{gordo_grass_par}
With $P$ and $\Ho_i$ defined
by (\ref{autoint_grass}) and (\ref{Hodge_grass}), the
following equality holds:

\begin{equation}\label{formulon_grass}
P=-\sum_{i=1}^{[{n\over2}]}{a_i\over
a_{i-1}+2a_i+a_{i+1}}\Ho_i
+\sum_{i=1}^{[{n+1\over2}]-1}{a_i
\begin{vmatrix}a_{i-1}+a_i&a_i+a_{i+1}\\
\alpha_i&\alpha_{i+1}\end{vmatrix}^2
\over
(a_{i-1}+2a_i+a_{i+1})(a_{i-1}+a_i)(a_i+a_{i+1})}
\end{equation}
\end{lma}

\begin{pf} 
We regard the expressions as polynomials in the indeterminates
$\alpha_1,\ldots,\alpha_{[{n+1\over2}]},
\lambda_1,\ldots,
\break 
\lambda_{[{n\over2}]}$
with coefficients in the quotient field of
$\C[a_1,\ldots,a_k]$. The only non-zero coefficients are
those of the monomials of the type 
$\alpha_i^2,\alpha_i\alpha_{i+1},\lambda_i$. Hence it
suffices to check that the corresponding coefficients in
the left and right hand sides coincide. This is trivial for
the coefficient of $\lambda_i$, which is $a_i$ in both sides.

For $\alpha_i\alpha_{i+1}$, its coefficient is zero in
the left-hand side of (\ref{formulon_grass}), while in the
right-hand side it is:
$$2{a_i\over
a_{i-1}+2a_i+a_{i+1}}-2{a_i\over
(a_{i-1}+2a_i+a_{i+1})(a_{i-1}+a_i)(a_i+a_{i+1})}(a_i+a_{i+1})(a_{i-1}+a_i)$$
which is also zero.

For $\alpha_i^2$, its coefficient in the left-hand side is $1$, while in the
right-hand side is: 
\begin{multline*}{a_{i-1}\over a_{i-2}+2a_{i-1}+a_i}
+{a_i\over
a_{i-1}+2a_i+a_{i+1}}+{a_{i-1}(a_{i-2}+a_{i-1})^2
\over(a_{i-2}+2a_{i-1}+a_i)(a_{i-2}+a_{i-1})(a_{i-1}+a_i)}\\
+{a_{i}(a_{i}+a_{i+1})^2\over
(a_{i-1}+2a_{i}+a_{i+1})(a_{i-1}+a_{i}) (a_{i}+a_{i+1})}
\end{multline*}
and a simple calculation shows that this is also $1$.

We leave to the reader to check that all these equalities remain
valid when $i$ takes the minimum and maximum values; a
special care should be taken when $n$ is even, because (as
observed at the end of Remark \ref{Si_irreducibles}), the
expression $a_{i-1}+2a_i+a_{i+1}$ appearing in both sums of
(\ref{formulon_grass}) must be regarded as
$(a_{i-1}+a_i)+(a_i+a_{i+1})$, thus just $a_{i-1}+a_i$ if
$i={n\over 2}$.
\end{pf}

\begin{cor}\label{rango-uno}
The matrix $M$ of (\ref{matrizota_grass}) has rank one.
\end{cor}

\begin{pf}
The left-hand side of (\ref{formulon_grass}) is zero by Lemma
\ref{lema_autoint_grass}, while the right-hand side is a sum of nonnegative
terms by (\ref{Hodge_grass}). Therefore we get equalities
$$\begin{vmatrix} a_{i-1}+a_i & \alpha_i+\alpha_{i+1}\\
\alpha_i & \lambda_i \end{vmatrix}+
\begin{vmatrix} a_i+a_{i+1} & \alpha_i+\alpha_{i+1}\\
\alpha_{i+1} & \lambda_i \end{vmatrix}=0$$
$$\begin{vmatrix}a_{i-1}+a_i&a_i+a_{i+1}\\
\alpha_i&\alpha_{i+1}\end{vmatrix}=0$$
(observe that for this it was crucial the hypothesis $a_i\ne0$ of Theorem
\ref{Barth-Larsen_grass}). The first equalities imply that the
columns of
$M$ containing a $\lambda_i$ depend on the other columns, while the
second inequalities show that the columns not containing any
$\lambda_i$ span a one-dimensional subspace.
\end{pf}

\section{Subvarieties of products of projective spaces}\label{prod-proj}

In this section we will apply our method to a product of two
projective spaces of the same dimension. The main result we
will prove is the following:

\begin{tma}\label{Barth-Larsen_PnxPn}
Let $X\subset \P^{n-1}\times\P^{n-1}$ be  smooth of dimension
greater than or equal to $n$ such that the two projection maps
$X\to\P^{n-1}$ are surjective. Then Pic$X$ is a free
abelian group of rank two, in which the pullbacks $H_1$ and
$H_2$ of the hyperplane classes of each $\P^{n-1}$ are
linearly independent.
\end{tma}

\begin{rem}{\rm
It is clear that the hypothesis of the surjectivity of the
projections $X\to\P^{n-1}$ are necessary. For instance, if
the image of the first projection is contained in a subvariety
$Z\subset\P^{n-1}$ with Picard group of rank at least two, then we can find
many smooth subvarieties $X\subset Z\times\P^{n-1}$ whose Picard group has
rank at least three (for example taking successive hyperplane sections of
the Segre embedding of $\P^{n-1}\times\P^{n-1}$). To give a concrete example,
just consider the embedding
$\P^1\times\P^1\times\P^1\times\P^1\subset\P^3\times\P^3$,
(in which we identify the first and last $\P^1\times\P^1$ with a smooth
quadric in $\P^3$ after the Segre embedding).

On the other hand, without any condition about
surjectivity, the results of \cite{S} show that Pic$X$ is
generated by $H_1$ and $H_2$ if $\dim X\ge\frac{3n-1}{2}$.

}\end{rem}

\begin{rem}\label{reduccion_PnxPn}{\rm
After the result of Debarre showing that $X$ is simply
connected \cite[Corollaire 2.4]{D}, as remarked in subsection
\ref{metodo_general} it is enough to assume $\dim X=n$ and
prove that every smooth divisor $D$ on $X$ restricted to
the intersection of $X$ with $n-2$ general hyperplanes
(i.e. divisors in the class $H_1+H_2$) is numerically
equivalent to a rational combination of
$H_1$ and $H_2$. Observe that $H_1,H_2$ are linearly
independent, since $H_1^n=H_2^n=0$ and $(H_1+H_2)^n>0$.
}\end{rem}

\subsection{General set-up} We fix now the general notation
that we will use throughout the section. Let $X$ be a smooth
$n$-dimensional subvariety of
$\P^{n-1}\times\P^{n-1}$ such that the two projections
$X\to\P^{n-1}$ are surjective. Let
$D$ be a smooth divisor of
$X$ and let $H_1,H_2$ be the pullback to $X$ of the
hyperplane classes of $\P^{n-1}$. 
We define $a_i$, $\alpha_i$, and $\lambda_i$ as the
intersection products in $X$ according to the following
table (sometimes it will be practical to set $a_0=a_n=0$ in order to get
compact expressions):
\begin{equation*}
\begin{tabular}{c|cccccccccc}
  &\small $H_1^{n-1}$ &\small $H_1^{n-2}H_2$ &\small $H_1^{n-3}H_2^2$ &\small
$\dots$ &\small
$H_2^{n-1}$ &\small $DH_1^{n-2}$ &\small $DH_1^{n-3}H_2$ &\small
$DH_1^{n-4}H_2^2$ &
$\dots$ &\small
$DH_2^{n-2}$ \\
\hline 
$\begin{matrix} H_1 \\ H_2 \\ D  \end{matrix}$ & 
$\begin{matrix} 0 \\ a_1 \\ \alpha_1  \end{matrix}$ & 
$\begin{matrix} a_1 \\ a_2 \\ \alpha_2  \end{matrix}$ & 
$\begin{matrix} a_2 \\ a_3 \\ \alpha_3  \end{matrix}$ &
$\begin{matrix} \dots \\ \dots \\ \dots  \end{matrix}$ & 
$\begin{matrix} a_{n-1} \\ 0 \\ \alpha_n  \end{matrix}$ & 
$\begin{matrix} \alpha_1 \\ \alpha_2 \\ \lambda_1  \end{matrix}$ & 
$\begin{matrix} \alpha_2 \\ \alpha_3 \\ \lambda_2  \end{matrix}$ & 
$\begin{matrix} \alpha_3 \\ \alpha_4 \\ \lambda_3  \end{matrix}$ & 
$\begin{matrix} \dots \\ \dots \\ \dots  \end{matrix}$ & 
$\begin{matrix} \alpha_{n-1} \\ \alpha_n \\ \lambda_{n-1}  \end{matrix}$ \end{tabular}
\end{equation*}

If $D$ depends numerically on $H_1$ and
$H_2$, it is clear that the intersection matrix 
\begin{equation}\label{matrizota}
M:=\begin{pmatrix}
0 & a_1 & a_2 & \ldots & a_{n-1} & \alpha_1 & \alpha_2 & \alpha_3 & \ldots &\alpha_{n-1} \\
a_1 & a_2 & a_3 & \ldots & 0 & \alpha_2 & \alpha_3 & \alpha_4 & \ldots &\alpha_{n} \\
\alpha_1 & \alpha_2 & \alpha_3 & \ldots & \alpha_{n} & \lambda_1 & \lambda_2 & \lambda_3 & \ldots &\lambda_{n-1}
\end{pmatrix}
\end{equation}
must have rank at least two. We claim that in fact if $M$ has rank two
then Theorem \ref{Barth-Larsen_PnxPn} follows. Indeed, observe first
that $\begin{vmatrix} 0 & a_{n-1} \\ a_1 & 0
\end{vmatrix} \ne 0$, since this is exactly the hypothesis that
the projections $X\to\P^{n-1}$ are surjective. Hence, if $M$ has rank two, it
follows that the third row is a linear combination of the first two.
This implies that there are $p,q\in\Q$ such that the product
of 
$$D':=D-pH_1-qH_2$$
with $H_1^{n-1},H_1^{n-2}H_1,\ldots,H_2^{n-1},
DH_1^{n-2},DH_1^{n-3}H_2,\ldots,DH_2^{n-2}$ is zero. If $S$
is the intersection of $X$ with $n-2$ general hyperplanes,
then it follows that
$$D'_{|S}(H_1+H_2)_{|S}=D'(H_1+H_2)^{n-1}=0$$ 
$${D'_{|S}}^2=D'(D-pH_1-qH_2)(H_1+H_2)^{n-2}=0.$$ 
Therefore, by Hodge index theorem $D'_{|S}$ is numerically
equivalent to zero. By Remark \ref{reduccion_PnxPn}, this
implies Theorem \ref{Barth-Larsen_PnxPn}.

We will thus devote the rest of the section to prove that $M$
has rank two (Corollary \ref{rango_matrizota_PnxPn}).

\subsection{First numerical properties}

In this subsection we prove the first relations between
the numerical invariants of a smooth subvariety
$X\subset\P^{n-1}\times\P^{n-1}$ of dimension $n$. 

Observe that the intersection matrix (\ref{matrizota}) shows
that the classes of $X$ and $D$ in the Chow ring of $\P^{n-1}\times\P^{n-1}$
are 
\begin{equation}\label{claseX}
[X]=a_{n-1}H_1^{n-2}+a_{n-2}H_1^{n-3}H_2+\ldots+a_1H_2^{n-2}
\end{equation}
\begin{equation}\label{claseD}
[D]=\alpha_nH_1^{n-1}+\alpha_{n-1}H_1^{n-2}H_2
+\ldots+\alpha_1H_2^{n-1}
\end{equation}
(we abuse momentarily the notation by indicating with
$H_1,H_2$ the hyperplane classes in each $\P^{n-1}$, instead
of denoting their pullback to $X$). This implies the
following:

\begin{lma}\label{autoint_lma_PnxPn}
The following equality holds
\begin{equation}\label{autoint_PnxPn}
P:=\alpha_1\alpha_{n}+\alpha_2\alpha_{n-1}+\ldots+\alpha_{n}\alpha_1-a_1\lambda_{n-1}-a_2\lambda_{n-2}-\ldots-a_{n-1}\lambda_1=0
\end{equation}
\end{lma}

\begin{pf}
This is (\ref{chern_general}) using  (\ref{claseX}), (\ref{claseD}) and the
intersection products given by the matrix $M$ of (\ref{matrizota}).
\end{pf}

A first important observation about our numerical invariants is that, in our
situation, all the numbers $a_1,\ldots,a_{n-1}$ are different from zero:

\begin{lma}\label{condic_numericas_PnxPn}
Let $X\subset \P^{n-1}\times\P^{n-1}$ be irreducible of
dimension $n$. Let $a_1,\ldots,a_{n-1}$ be defined as above.
Suppose that $a_1$, $a_{n-1}$ are different from zero. Then
all $a_i$ are strictly positive.
\end{lma}

\begin{pf} 
See \cite[Th\'eor\`eme 1.3]{D}.
\end{pf}

\begin{df} For $i=1,\ldots,n-1$, we will denote $S_i$ to the
surface obtained intersecting $X$ with a general
$\P^i\times\P^{n-i}$ (hence the class of $S_i$ in $X$ is
$H_1^{n-i-1}H_2^{i-1}$). Since both natural projections
restricted to $X$ are surjective, it follows from
\cite[Th\'eor\`eme 1.3]{D}  that $S_i$ is a smooth
irreducible surface.
\end{df}

\begin{rem}\label{relaciones_a_alpha_PnxPn}
{\rm
The Hodge index theorem applied to the divisors $H_1,H_2$ 
restricted to $S_i$ implies the inequality:
$$\begin{vmatrix} a_{i-1} & a_i \\ a_i & a_{i+1}
\end{vmatrix}=\begin{vmatrix}{H_1}_{|S_i}^2 &
{H_1}_{|S_i}{H_2}_{|S_i}\\ {H_1}_{|S_i}{H_2}_{|S_i}&
{H_2}_{|S_i}^2\end{vmatrix}\le 0$$ 
with equality if and only
if $H_1$ and $H_2$ are numerically dependent on $S_i$.
Lemma \ref{condic_numericas_PnxPn} allows us to write all
these inequalities as:
\begin{equation}\label{cadena_ai_PnxPn}
{a_1\over a_2}\le {a_2\over a_3}\le\ldots\le{a_{n-2} \over
a_{n-1}}.
\end{equation}

Thus it follows that, for all $i<j$, we have
$$\begin{vmatrix} a_{i-1} & a_j \\ a_i & a_{j+1}
\end{vmatrix} \le 0$$ with equality if and only if
$$\begin{vmatrix} a_{i-1} & a_i \\ a_i & a_{i+1}
\end{vmatrix}=\begin{vmatrix} a_i & a_{i+1} \\ a_{i+1} &
a_{i+2}\end{vmatrix}=\ldots=\begin{vmatrix} a_{j-1} & a_j \\
a_j & a_{j+1} \end{vmatrix} = 0$$ i.e.
$${a_{i-1}\over a_i}={a_i\over a_{i+1}}=\ldots={a_j
\over a_{j+1}}.$$ 
We suspect that the inequalities in (\ref{cadena_ai_PnxPn})
are never equalities. However we are not able to prove it,
so that our proofs will become sometimes more complicated when
considering the cases in which some equalities hold.
}\end{rem}

\begin{rem}\label{Hodge-PnxPn}{\rm
We apply now the Hodge index theorem to the restriction to
$S_i$ of the divisors $H_1,H_2,D$ and get

\begin{equation}\label{Hodge-grande}
r_{i,i+1,n+i}:=
\begin{vmatrix}
a_{i-1} & a_{i} & \alpha_{i} \\
a_{i} & a_{i+1} & \alpha_{i+1} \\
\alpha_{i} & \alpha_{i+1} & \lambda_i
\end{vmatrix}=
\begin{vmatrix}
H_{1|S_i}^2 & H_{1|S_i}H_{2|S_i} & H_{1|S_i}D_{|S_i}\\
H_{1|S_i}H_{2|S_i} & H_{2|S_i}^2 & H_{2|S_i}D_{|S_i}\\
H_{1|S_i}D_{|S_i} & H_{2|S_i}D_{|S_i} & D_{|S_i}^2
\end{vmatrix}\ge 0
\end{equation}
(in general $r_{i,j,k}$ will denote the minor of the columns
$i,j,k$ of the matrix $M$) with equality if and only if $H_1$,
$H_2$ and $D$ are numerically dependent when restricted to $S_i$. Hence the
vanishing of $r_{i,i+1,n+i}$ would immediately prove that the
restriction of $D$ to $S_i$ is numerically a combination of
$H_1$ and $H_2$, unless the restrictions to
$S_i$ of $H_1$ and $H_2$ are numerically dependent. 

In this latter case (which occurs if and only if ${a_{i-1}\over
a_i}={a_i\over a_{i+1}}$, by Remark
\ref{relaciones_a_alpha_PnxPn}) the restriction of $D$ to
$S_i$ is numerically a combination of $H_1$ and $H_2$ if and
only if it is just a multiple of $H_1$. By the Hodge index
theorem, this is equivalent to the inequality:
\begin{equation}\label{Hodge-reducida}
\tilde r_i:=\begin{vmatrix}
a_{i-1} & \alpha_{i} \\
\alpha_{i} & \lambda_i
\end{vmatrix}=\begin{vmatrix}
H_{1|S_i}^2 & H_{1|S_i}D_{|S_i}\\
H_{1|S_i}D_{|S_i} & D_{|S_i}^2
\end{vmatrix}\le 0
\end{equation}
to be an equality. Moreover, it easily follows that
$$\mathrm{rk}\begin{pmatrix}
a_{i-1} & a_{i} & \alpha_{i} \\
a_{i} & a_{i+1} & \alpha_{i+1}
\end{pmatrix}=1$$
since this is the multiplication matrix in $S_i$ of the
divisors $H_1,H_2$ (which are numerically dependent) with
$H_1,H_2,D$. In particular 
\begin{equation}\label{primera-sustitucion-alphas}
 \alpha_{i+1}={a_{i+1}\over a_{i}}\alpha_i
\end{equation}
}\end{rem}

\begin{rem} {\rm The above remark shows that we should
check how many equalities we have in the chain
(\ref{cadena_ai_PnxPn}). A practical way to control this is to
define the function:
\begin{equation}\label{sigma_PnxPn}
\begin{matrix}
\sigma : & \{ 2,\ldots,n\} & \longrightarrow & \{
2,\ldots,n\}\\
 & i & \mapsto &
\begin{matrix}\mathrm{unique}\ l\ \mathrm{such\ that}\\
{a_{l-2}\over a_{l-1}}<{a_{l-1}\over a_l}=\ldots={a_{i-2}\over
a_{i-1}}={a_{i-1}\over a_i}\end{matrix}
\end{matrix}
\end{equation}
In this way, the equality ${a_{i-1}\over
a_i}={a_{j-1}\over a_j}$ is equivalent to the equality $\sigma(i)=\sigma(j)$.
In other words, the equalities in (\ref{cadena_ai_PnxPn})
coincide with the equalities in 
$$\sigma(2)\le\sigma(3)\le\ldots\le\sigma(n).$$
Observe that saying that there are no equalities in (\ref{cadena_ai_PnxPn})
(as we said, this is what we suspect) is equivalent to saying that
$\sigma$ is the identity map.

We recall from Remark \ref{relaciones_a_alpha_PnxPn} that
${a_{\sigma(i)-2}\over a_{\sigma(i)-1}}<{a_{\sigma(i)-1}\over
a_{\sigma(i)}}=\ldots={a_{i-2}\over a_{i-1}}={a_{i-1}\over
a_i}$ is equivalent to 
$\begin{vmatrix} a_{\sigma(i)-1} & a_{i-1} \\ 
a_{\sigma(i)} & a_i \end{vmatrix}=0,\ \ \ 
\begin{vmatrix} a_{\sigma(i)-2} & a_{i-1} \\ 
a_{\sigma(i)-1} & a_i \end{vmatrix} < 0$. On the other
hand, iterating (\ref{primera-sustitucion-alphas}) we get 
\begin{equation}\label{sustitucion-alphas}
\alpha_{i}={a_{i}\over a_{\sigma(i)}}\alpha_{\sigma(i)}
\end{equation}
}\end{rem}

The idea of the above function $\sigma$ is that it allows to
determine ``virtual'' equalities in the chain
(\ref{cadena_ai_PnxPn}), and then perform the substitutions
given in (\ref{sustitucion-alphas}). This will be very useful
in the next section, when $a_1,\ldots,a_{n-1}$ will be
indeterminates rather than concrete values. This motivates
the following definition, which indicates when a function
$\sigma$ is obtained as in the above remark:

\begin{df}
We call {\it partition map} to a function $\sigma: \{
2,\ldots,n\} \to \{ 2,\ldots,n\}$ 
satisfying:
\begin{enumerate}
\item $\sigma(i)\le i$ for all $i$  
\item for all $j \in \{\sigma(i),\ldots,i\}$, $\sigma(j)=\sigma(i)$. 
\end{enumerate}
\end{df}

\subsection{Polynomial equalities}

Throughout this subsection, we will consider $a_i$,
$\alpha_i$, $\lambda_i$ as indeterminates, and we will
consider the polynomial ring
$\Q(a_1,\ldots,a_{n-1})[\alpha_1,\ldots,
\alpha_n,\lambda_1,\ldots,\lambda_{n-1}]$. In this context,
the elements of $\Q(a_1,\ldots,a_{n-1})$ will be regarded
as coefficients (our convention $a_0=a_n=0$ will be useful again for writing
the expressions in a compact way). We define the polynomials (compare with
(\ref{autoint_PnxPn}), (\ref{Hodge-grande}) and (\ref{Hodge-reducida})): 
$$P:=\alpha_1\alpha_{n}+\alpha_2\alpha_{n-1}+
\ldots+\alpha_{n}\alpha_1-a_1\lambda_{n-1}
-a_2\lambda_{n-2}-\ldots-a_{n-1}\lambda_1$$
$$r_{i,i+1,n+i}:=
\begin{vmatrix}
a_{i-1} & a_{i} & \alpha_{i} \\
a_{i} & a_{i+1} & \alpha_{i+1} \\
\alpha_{i} & \alpha_{i+1} & \lambda_i
\end{vmatrix}$$
$$\tilde r_i:=\begin{vmatrix}
a_{i-1} & \alpha_{i} \\
\alpha_{i} & \lambda_i
\end{vmatrix}$$
$$r_{i,i+1,n-i+1}:=
\begin{vmatrix}
a_{i-1} & a_{i} & a_{n-i} \\
a_{i} & a_{i+1} & a_{n-i+1} \\
\alpha_{i} & \alpha_{i+1} & \alpha_{n-i+1}
\end{vmatrix}$$
and the coefficients
\begin{equation}\label{def_c_i}
c_i:={a_{n-i}\over
\begin{vmatrix} a_{i-1} & a_i
\\ a_i & a_{i+1} \end{vmatrix}}
\end{equation}
\begin{equation}\label{def_d_i}
d_i:={a_ia_{n-i}\over
\begin{vmatrix} a_{i-1} & a_i
\\ a_i & a_{i+1} \end{vmatrix} \begin{vmatrix}
a_{i-1} & a_{n-i} \\ a_{i} & a_{n-i+1} \end{vmatrix}
\begin{vmatrix} a_{i} & a_{n-i}
\\ a_{i+1} & a_{n-i+1} \end{vmatrix} }
\end{equation}
Observe that, if $i=\lceil {n\over2}\rceil$ then
$r_{i,i+1,n-i+1}=0$ and $d_i$ is undefined. We will take both of them equal to zero.

\begin{lma}\label{caso0_PnxPn} The following equality holds in
$\Q(a_1,\ldots,a_{n-1})[\alpha_1,\ldots,\alpha_n,\lambda_1,\ldots,
\lambda_{n-1}]$:
\begin{equation}\label{formulon2_PnxPn}
P=-\sum_{i=1}^{n-1}c_ir_{i,i+1,n+i}
-\sum_{i=1}^{n-1}d_ir_{i,i+1,n-i+1}^2
\end{equation}
\end{lma}

\begin{pf}
We have to compare the coefficients of the different
monomials. The only monomials appearing in
(\ref{formulon2_PnxPn}) are of the type
$\lambda_i,\alpha_i^2,\alpha_i\alpha_{i+1},\alpha_i\alpha_{n-i+1}$
and $\alpha_i\alpha_{n-i+2}$. It is immediate that the
coefficient of $\lambda_i$ is $-a_{n-i}$ in both the left
and right hand sides of (\ref{formulon2_PnxPn}).

For $\alpha_i^2$ (we do not consider the case
$i= 1,\lceil{n\over2}\rceil,\lceil{n\over2}\rceil+1,n$, which overlaps with other cases),
we have zero on the left-hand side, while the coefficient on
the right-hand side comes from
$r_{i-1,i,n+i-1}$, $r_{i,i+1,n+i}$, $r_{i-1,i,n-i+2}^2$,
$r_{n-i+1,n-i+2,i}^2$ and $r_{i,i+1,n-i+1}^2$ and it is
$$c_{i-1}a_{i-2}+c_ia_{i+1}
-d_{i-1}\begin{vmatrix} a_{i-2} & a_{n-i+1} \\ 
a_{i-1} & a_{n-i+2} \end{vmatrix}^2-
d_{i}\begin{vmatrix} a_{i} & a_{n-i} \\ 
a_{i+1} & a_{n-i+1}\end{vmatrix}^2-
d_{n-i+1}\begin{vmatrix} a_{n-i} & a_{n-i+1} \\
a_{n-i+1}  & a_{n-i+2} \end{vmatrix}^2$$
which is also zero.

For $\alpha_i\alpha_{i+1}$, (we assume now $i\ne \lceil{n\over2}\rceil,\lceil {n\over2}\rceil+1$ in order to avoid
overlaps with other cases) we have again zero on the
left-hand side, while the coefficient on the right-hand
side comes from
$r_{i,i+1,n+i}$ and $r_{i,i+1,n-i+1}^2$, and it is
$$-2c_ia_i+2d_i\begin{vmatrix} a_{i-1} & a_{n-i} \\ a_i &
a_{n-i+1}
\end{vmatrix}\begin{vmatrix} a_i & a_{n-i} \\ 
a_{i+1} & a_{n-i+1} \end{vmatrix}$$
which is obviously zero.

For $\alpha_i\alpha_{n-i+1}$ (assuming now
$i\ne 1,\lceil{n\over2}\rceil,\lceil{n\over2}\rceil+1$), the coefficient on the
left-hand side is $2$, while on the right-hand side, it
comes from $r_{i,i+1,n-i+1}^2$ and $r_{n-i+1,n-i+2,i}^2$,
and it is
$$-2d_i\begin{vmatrix} a_{i} & a_{n-i} \\ a_{i+1} & a_{n-i+1}
\end{vmatrix}\begin{vmatrix} a_{i-1} & a_{i} \\ 
a_{i} & a_{i+1} \end{vmatrix}
+2d_{n-i+1}\begin{vmatrix} a_{n-i} & a_{n-i+1} \\
a_{n-i+1} & a_{n-i+2} \end{vmatrix}
\begin{vmatrix} a_{i-1} & a_{n-i+1} \\ 
a_{i} & a_{n-i+2} \end{vmatrix}$$ 
which is also $2$.

For $\alpha_i\alpha_{n-i+2}$ (again with
$i\ne\lceil{n\over2}\rceil,\lceil{n\over2}\rceil+1$), its coefficient is zero in the
left-hand side, while it appears only in the right-hand
side in $r_{i-1,i,n-i+2}^2$ and $r_{n-i+1,n-i+2,i}^2$, with
coefficient: 
$$2d_{i-1}\begin{vmatrix} a_{i-2} & a_{n-i+1} \\ 
a_{i-1} & a_{n-i+2} \end{vmatrix}
\begin{vmatrix} a_{i-2} & a_{i-1} \\ 
a_{i-1} & a_{i} \end{vmatrix}
+2d_{n-i+1}\begin{vmatrix} a_{n-i} & a_{n-i+1} \\ 
a_{n-i+1} & a_{n-i+2} \end{vmatrix}
\begin{vmatrix} a_{i-1} & a_{n-i} \\ 
a_{i} & a_{n-i+1} \end{vmatrix}$$
which is zero.

We are left with the coefficients of the cases we have excluded, i.e. the ones of the monomials $\alpha_1^2,\alpha_1\alpha_n, \alpha_k^2,\alpha_k\alpha_{k+1}, \alpha_k\alpha_{k+2},\alpha_{k+1}^2,\alpha_{k+1}\alpha_{k+2}$
and $\alpha_n^2$, with $k=\lceil{n\over2}\rceil$. We omit the proof, which is as before a simple computation, but depending now on
whether $n=2k$ or $n=2k-1$ for the cases involving $k$ in an index.
\end{pf}

Thinking of our application to subvarieties of
$\P^{n-1}\times\P^{n-1}$, Lemma \ref{caso0_PnxPn} could be
the analogue of Lemma \ref{gordo_grass_par}, but a priori
the denominators of some coefficient $c_i,d_i$ could
vanish, depending on whether we have equalities in the chain
(\ref{cadena_ai_PnxPn}). If this happens, the idea would be to use, when
necessary, the inequalities (\ref{Hodge-reducida}) instead of
(\ref{Hodge-grande}) and perform the substitutions
(\ref{sustitucion-alphas}). In order to keep all these substitutions in the
framework of polynomials, we fix a partition map
$\sigma$ for the rest of the subsection.

The first case in which we would be in trouble is when the
denominator of $c_i$ vanishes (i.e. when
$\sigma(i)=\sigma(i+1)$, in terms of the partition map). As we
observed in Remark \ref{Hodge-PnxPn}, in this case we
should replace $r_{i,i+1,n+i}$ with $\tilde r_i$ (with an
appropriate coefficient to fit in (\ref{formulon2_PnxPn})).
On the other hand, the denominator of $d_i$ would also
create problems, so it looks convenient to remove the term
$d_ir_{i,i+1,n-i+1}^2$ from (\ref{formulon2_PnxPn}). The
precise way of performing these changes is given by the
following:

\begin{lma}\label{cambio_de_caso_PnxPn} Suppose 
$\sigma(i)=\sigma(i+1)$ for some $i=1,\ldots,n-1$ and
set:
\begin{equation*}
l_i:={
\begin{matrix}a_{n-i}\end{matrix}
\begin{vmatrix}  a_{i} & a_{n-i+1}\\
\alpha_i & \alpha_{n-i+1} \end{vmatrix}^2 
\begin{vmatrix} a_{i-1} & a_{i} \\ 
a_i & a_{i+1} \end{vmatrix} \over 
\begin{matrix}a_{i}\end{matrix}
\begin{vmatrix} a_{i-1} & a_{n-i} \\ 
a_i & a_{n-i+1} \end{vmatrix}
\begin{vmatrix} a_{i} & a_{n-i} \\ 
a_{i+1} & a_{n-i+1} \end{vmatrix}} 
\end{equation*}
\begin{equation*}
m_i:={a_{n-i}\begin{vmatrix} a_i & a_{i+1} \\ \alpha_i & \alpha_{i+1}
\end{vmatrix} (a_ia_{n-i+1}\alpha_{i+1}
-2a_ia_{i+1}\alpha_{n-i+1}+a_{i+1}a_{n-i+1}\alpha_i) \over a_ia_{i+1}
\begin{vmatrix} a_i & a_{n-i} \\ a_{i+1} & a_{n-i+1} \end{vmatrix}}
\end{equation*}
Then, 
$$
-c_ir_{i,i+1,n+i}-d_ir_{i,i+1,n-i+1}^2=-{a_{n-i}\over
a_{i+1}}\tilde r_i-l_i-m_i
$$
As a consequence, 
\begin{multline}\label{formulon_gral_PnxPn}
P = -\sum_{\sigma(i)<\sigma(i+1)}
c_ir_{i,i+1,n+i}- \sum_{\sigma(i)=\sigma (i+1)}
{a_{n-i}\over a_{i+1}} \tilde r_i-\\ 
-\sum_{\sigma(i)<\sigma(i+1)}
d_ir_{i,i+1,n-i+1}^2
-\sum_{\sigma(i)=\sigma (i+1)}l_i
-\sum_{\sigma(i)=\sigma (i+1)}m_i
\end{multline}

\end{lma}

\begin{pf} The first part is a straightforward computation,
while the second one follows immediately from the first one
and (\ref{formulon2_PnxPn}).
\end{pf}

Observe that we can still have problems with the
denominator of $l_i$ and $m_i$ if $\sigma(n-i+1)$ equals
$\sigma(i)=\sigma(i+1)$. In the same way, even if
$\sigma(i)<\sigma(i+1)$ we can have problems with the
denominator of $d_i$ if $\sigma(n-i+1)$ equals
$\sigma(i)$ or $\sigma(i+1)$. Let us see that we can
avoid these problems by performing the substitutions
suggested by (\ref{sustitucion-alphas}). We thus make the
following: 

\begin{df}
Given any polynomial
$Q\in\Q(a_1,\ldots,a_{n-1})[\alpha_1,\ldots,\alpha_n,\lambda_1,\ldots,\lambda_{n-1}]$,
we define a new polynomial (depending on few variables): 
$$Q^\sigma:=Q(a_1,\ldots,a_{n-1},\alpha_{1},{a_{2}\over
a_{\sigma(2)}}\alpha_{\sigma(2)},\ldots,{a_{n-1}\over
a_{\sigma(n-1)}}\alpha_{\sigma(n-1)},\alpha_n,
\lambda_1,\ldots,\lambda_{n-1})$$
\end{df}

We see now how these substitutions affect the problematic
terms in (\ref{formulon_gral_PnxPn}).

\begin{lma}\label{incordiones_PnxPn} The following
identities hold:
\item{(i)} If $\sigma(i)=\sigma(i+1)$, then $m_i^\sigma=0$.
\item{(ii)} If $\sigma(i)=\sigma(i+1)=\sigma(n-i+1)$, then
$l_i^\sigma=0$.
\item{(iii)} If $\sigma(i)<\sigma(i+1)=\sigma(n-i+1)$, then
$r_{i,i+1,n-i+1}^\sigma={-a_i\over
a_{\sigma(i)}a_{i+1}}
\begin{vmatrix} a_{i} & a_{n-i} \\
a_{i+1} & a_{n-i+1} \end{vmatrix}
\begin{vmatrix} a_{\sigma(i)} & a_{i+1} \\
\alpha_{\sigma(i)} & \alpha_{i+1} \end{vmatrix}$
and hence 
\begin{equation}\label{def_g_i}
g_i:=d_i(r_{i,i+1,n-i+1}^\sigma)^2=
{a_i^3a_{n-i}
\begin{vmatrix} a_i & a_{n-i} \\ 
a_{i+1} & a_{n-i+1}\end{vmatrix}
\begin{vmatrix} a_{\sigma(i)} & a_{i+1} \\
\alpha_{\sigma(i)} & \alpha_{i+1} \end{vmatrix}^2
\over a_{\sigma(i)}^2a_{i+1}^2
\begin{vmatrix} a_{i-1} & a_{n-i} \\ 
a_{i} & a_{n-i+1} \end{vmatrix}
\begin{vmatrix} a_{i-1} & a_i \\ 
a_i & a_{i+1} \end{vmatrix}}
\end{equation}
\item{(iv)} If $\sigma(n-i+1)=\sigma(i)<\sigma(i+1)$, then
$r_{i,i+1,n-i+1}^\sigma={-1\over \begin{matrix}a_{\sigma(i)}\end{matrix}}
\begin{vmatrix} a_{i-1} & a_{n-i} \\ 
a_i & a_{n-i+1}\end{vmatrix}
\begin{vmatrix} a_{\sigma(i)} & a_{i+1} \\
\alpha_{\sigma(i)} & \alpha_{i+1} \end{vmatrix}$ 
and hence 
\begin{equation}\label{def_h_i}
h_i:=d_i(r_{i,i+1,n-i+1}^\sigma)^2=
{a_ia_{n-i}\begin{vmatrix} a_{i-1} & a_{n-i} \\ a_{i} &
a_{n-i+1}\end{vmatrix}\begin{vmatrix} a_{\sigma(i)} & a_{i+1} \\
\alpha_{\sigma(i)} & \alpha_{i+1} \end{vmatrix}^2\over
a_{\sigma(i)}^2\begin{vmatrix} a_{i} & a_{n-i} \\ a_{i+1} &
a_{n-i+1} \end{vmatrix}\begin{vmatrix} a_{i-1} & a_i \\ a_i
& a_{i+1} \end{vmatrix}}
\end{equation}
\end{lma}

\begin{pf}
This is just a straightforward computation.
\end{pf}

\begin{cor}\label{Psigma}
 With the above notations, 
\begin{multline}\label{sigma_formulon_gral_PnxPn}
P^\sigma = -\sum_{\sigma(i)<\sigma(i+1)}
c_ir_{i,i+1,n+i}^\sigma - \sum_{\sigma(i)=\sigma (i+1)}
{a_{n-i}\over a_{i+1}} \tilde r_i^\sigma-\\
-\sum_{\sigma(n-i+1)\ne\sigma(i)<\sigma(i+1)\ne
\sigma(n-i+1)} d_i(r_{i,i+1,n-i+1}^\sigma)^2
-\sum_{\sigma(i)<\sigma(i+1)=
\sigma(n-i+1)} g_i
 -\\
-\sum_{\sigma(n-i+1)=\sigma(i)<\sigma(i+1)}
h_i-\sum_{\sigma(i)=\sigma
(i+1)\ne\sigma(n-i+1)}l_i^\sigma
\end{multline}
\end{cor}

\begin{pf}
Just apply Lemma \ref{incordiones_PnxPn} to the equality
(\ref{formulon_gral_PnxPn}) of Lemma
\ref{cambio_de_caso_PnxPn}.
\end{pf}

\subsection{End of the proof}

After the previous technical subsection we come back to the
geometrical situation of a subvariety
$X\subset\P^{n-1}\times\P^{n-1}$ with the usual invariants. We can now state
and prove the analogue of Lemma \ref{gordo_grass_par}. 

\begin{lma}\label{autoint-cambiada}
The following equality holds:
\begin{multline}\label{igualdad_de_cargarse_Hodges_PnxPn}
0 = -\sum_{{a_{i-1}\over a_i}<{a_i\over a_{i+1}}}
{\begin{matrix}a_{n-i}\end{matrix}\over \begin{vmatrix} a_{i-1}
& a_i
\\ a_i & a_{i+1}
\end{vmatrix}} \begin{vmatrix}
a_{i-1} & a_{i} & \alpha_{i} \\
a_{i} & a_{i+1} & \alpha_{i+1} \\
\alpha_{i} & \alpha_{i+1} & \lambda_i
\end{vmatrix} 
-\sum_{{a_{i-1}\over a_i}={a_i\over
a_{i+1}}} {a_{n-i}\over a_{i+1}} \begin{vmatrix} a_{i+1} &
\alpha_{i+1} \\ \alpha_{i+1} & \lambda_{i}
\end{vmatrix}-\\
-\sum_{\frac{a_{n-i}}{a_{n-i+1}}\ne\frac{a_{i-1}}{a_i}
<\frac{a_i}{a_{i+1}}\ne\frac{a_{n-i}}{a_{n-i+1}}}
{\begin{matrix}  a_ia_{n-i}\end{matrix}\over 
\begin{vmatrix}a_{i-1} & a_i 
\\ a_i & a_{i+1} \end{vmatrix} 
\begin{vmatrix} a_{i} & a_{n-i} \\ 
a_{i+1} & a_{n-i+1} \end{vmatrix} 
\begin{vmatrix} a_{i-1} & a_{n-i} \\ 
a_{i} & a_{n-i+1} \end{vmatrix}}
\begin{vmatrix}
a_{i-1} & a_{i} & a_{n-i} \\
a_{i} & a_{i+1} & a_{n-i+1} \\
\alpha_{i} & \alpha_{i+1} & \alpha_{n-i+1}
\end{vmatrix}^2
\end{multline}
and therefore
\begin{equation}\label{Hodge1_PnxPn_cero}
\begin{matrix}
\begin{vmatrix}
a_{i-1} & a_{i} & \alpha_{i} \\
a_{i} & a_{i+1} & \alpha_{i+1} \\
\alpha_{i} & \alpha_{i+1} & \lambda_i
\end{vmatrix}=0
& {\mathrm{ if }}\ {a_{i-1}\over a_i}<{a_i\over a_{i+1}}
\end{matrix}
\end{equation}
\begin{equation}\label{Hodge2_PnxPn_cero}
\begin{matrix}
\begin{vmatrix} a_{i+1} & \alpha_{i+1}\\ \alpha_{i+1} &
\lambda_i \end{vmatrix}=0  & {\mathrm{ if }}\ {a_{i-1}\over
a_i}={a_i\over a_{i+1}}
\end{matrix}
\end{equation}
\begin{equation}\label{menor_cero}
\begin{matrix}
\begin{vmatrix}
a_{i-1} & a_{i} & a_{n-i} \\
a_{i} & a_{i+1} & a_{n-i+1} \\
\alpha_{i} & \alpha_{i+1} & \alpha_{n-i+1}
\end{vmatrix}=0& {\mathrm{ if }}\
\frac{a_{n-i}}{a_{n-i+1}}\ne\frac{a_{i-1}}{a_i}<\frac{a_i}{a_{i+1}}\ne\frac{a_{n-i}}{a_{n-i+1}}
\end{matrix}
\end{equation}

\end{lma}

\begin{pf} Let $\sigma$ be the partition map defined by (\ref{sigma_PnxPn}).
Then (\ref{igualdad_de_cargarse_Hodges_PnxPn}) is an easy consequence of the
equality  (\ref{sigma_formulon_gral_PnxPn}) in Corollary \ref{Psigma}.
Indeed, observe first that, by (\ref{sustitucion-alphas}), any polynomial
$Q\in\Q(a_1,\ldots,a_{n-1})[\alpha_1,\ldots,\alpha_n,\lambda_1,\ldots,\lambda_{n-1}]$,
takes the same value as $Q^\sigma$ when applied to the invariants of $X$,
provided that all denominators are different from zero. In particular,
$P^\sigma$ vanishes thanks to Lemma \ref{autoint_lma_PnxPn}. We thus apply
Corollary \ref{Psigma} and analyze each of the summands of
(\ref{sigma_formulon_gral_PnxPn}): 

It is clear that the first three summands of
(\ref{sigma_formulon_gral_PnxPn}) coincide with the summands of
(\ref{igualdad_de_cargarse_Hodges_PnxPn}) and that their
denominators do not vanish (recall that this is equivalent to
say $\sigma(i)<\sigma(i+1)$). 
For the other summands of (\ref{sigma_formulon_gral_PnxPn}), if
$\frac{a_{i-1}}{a_i}<\frac{a_i}{a_{i+1}}=\frac{a_{n-i}}{a_{n-i+1}}$
or
$\frac{a_{n-i}}{a_{n-i+1}}=\frac{a_{i-1}}{a_i}<\frac{a_i}{a_{i+1}}$,
then $d_i(r_{i,i+1,n-i+1}^\sigma)^2$ takes the value $g_i$
or $h_i$ respectively of (\ref{def_g_i}) and
(\ref{def_h_i}), whose denominators do not vanish, but
their numerators contain respectively the
factor $\begin{vmatrix} a_{i-1} & a_{n-i} \\  a_i &
a_{n-i+1}\end{vmatrix}$ and
$\begin{vmatrix} a_{i} & a_{n-i} \\ a_{i+1} & a_{n-i+1}
\end{vmatrix}$, so that they vanish; and finally, if
$\frac{a_{i-1}}{a_i}=\frac{a_i}{a_{i+1}}\ne\frac{a_{n-i}}{a_{n-i+1}}$,
then the denominator of $l_i$, as defined in Lemma
\ref{cambio_de_caso_PnxPn} does not vanish, while the
numerator contains the factor $\begin{vmatrix} a_{i-1} &
a_i \\  a_i & a_{i+1}\end{vmatrix}$, hence it vanishes.

Observe now that all the summands of
(\ref{igualdad_de_cargarse_Hodges_PnxPn}) are nonnegative by
Remark \ref{Hodge-PnxPn}, which implies the equalities
(\ref{Hodge1_PnxPn_cero}), (\ref{Hodge2_PnxPn_cero}) and
(\ref{menor_cero}) (again, it was crucial to know that all $a_i$ are not
zero).
\end{pf}

We are now ready to prove that the matrix $M$ of
(\ref{matrizota}) has rank two. We will do it in several steps.

\begin{cor}\label{rango_dos_peque} 
For each $i=1,\ldots,n$, $i\ne\lceil {n\over 2}\rceil$, define the following submatrix of $M$:
$$M_i':=\left\{\begin{matrix}
\begin{pmatrix}
a_{i-1} & a_{i} & a_{n-i} \\
a_{i} & a_{i+1} & a_{n-i+1} \\
\alpha_{i} & \alpha_{i+1} & \alpha_{n-i+1}
\end{pmatrix}
 \ \ \ {\rm if}\ i<n-i\ \ \ \ \  \\ \\
\begin{pmatrix}
a_{n-i} & a_{i-1} & a_{i}  \\
a_{n-i+1} & a_{i} & a_{i+1}  \\
\alpha_{n-i+1} & \alpha_{i} & \alpha_{i+1} 
\end{pmatrix}
 \ \ \ {\rm if}\ n-i<i-1
\end{matrix}\right.
$$
Then the middle column of $M'_i$ is a linear combination of the
two others.
\end{cor}

\begin{pf}
This is immediate if
$\frac{a_{n-i}}{a_{n-i+1}}\ne\frac{a_{i-1}}{a_i}<\frac{a_i}{a_{i+1}}\ne\frac{a_{n-i}}{a_{n-i+1}}$,
since in this case any two columns of $M'_i$ are linearly
independent, and (\ref{menor_cero}) implies that any of the
columns of $M'_i$ is a linear combination of the two others. It
is also immediate if ${a_{i-1}\over a_i}={a_i\over a_{i+1}}$,
since in this case $\begin{pmatrix}a_{i-1}\\ a_i\\
\alpha_i\end{pmatrix}$ and $\begin{pmatrix}a_i\\ a_{i+1}\\
\alpha_{i+1}\end{pmatrix}$ are multiple of each
other (by the second part of Remark \ref{Hodge-PnxPn}) and one of them is the
middle column. We are thus left with the cases
$\frac{a_{i-1}}{a_i}<\frac{a_i}{a_{i+1}}=\frac{a_{n-i}}{a_{n-i+1}}$
and
$\frac{a_{n-i}}{a_{n-i+1}}=\frac{a_{i-1}}{a_i}<\frac{a_i}{a_{i+1}}$,
which are proved in the same way using again Remark
\ref{Hodge-PnxPn}. 
\end{pf}

\begin{cor}\label{rango_dos_media}
The submatrix 
$$M':=\begin{pmatrix}
0 & a_1 & a_2 & \ldots & a_{n-2} & a_{n-1} \\ 
a_1 & a_2 & a_3 & \ldots & a_{n-1} & 0 \\
\alpha_1 & \alpha_2 & \alpha_3 & \ldots & \alpha_{n-1} &
\alpha_n
\end{pmatrix}$$
of $M$ has rank two.
\end{cor}

\begin{pf}
It is enough to prove for each $i=1,\ldots,n$ that the $i$-th
column of $M'$ is a linear combination of the first and last
columns (which in turn are linearly independent). We will prove it
by induction on $d(i):=\min\{i-1,n-i\}$ (i.e. the distance of
$i$ to the border of the interval $[1,n]$). Of course, the
statement is trivial for $d(i)=0$ (i.e. $i=1$ or $i=n$). 

Assume first $i<{n+2\over2}$ (hence $d(i)=i-1$). Then Corollary
\ref{rango_dos_peque} implies that the $i$-th column is a
linear combination of the $(i-1)$-th and $(n-i+2)$-th columns.
By induction hypothesis (observe that $d(i-1)=d(n-i+2)=i-2$),
these two columns are in turn linear combinations of the first
and last columns of $M$, which completes the proof of this case.

If instead $i>{n+1\over2}$ (hence $d(i)=n-i$), Corollary
\ref{rango_dos_peque} implies now that the $i$-th column is a
linear combination of the $(n-i+1)$-th and $(i+1)$-th columns. The
$(i+1)$-th column is, again by induction hypothesis, a linear combination of
the first and last columns of $M$. For the $(n-i+1)$-th column we observe
that
$d(n-i+1)=n-i=d(i)$ and that $n-i+1<{n+2\over2}$. Hence we are
in the situation of the previous case, for which we already
proved the inductive step.
\end{pf}

\begin{cor}\label{rango_matrizota_PnxPn}
The matrix $M$ of (\ref{matrizota}) has rank two.
\end{cor}

\begin{pf}
By Corollary \ref{rango_dos_media}, it is enough to prove that
any column of $M$ containing a $\lambda_i$ is a linear
combination of the columns of $M'$. To this purpose, for
$i=1,\ldots, n-1$ we consider the submatrix:
$$M_i:=\begin{pmatrix}a_{i-1} & a_i & \alpha_i \\ 
a_i & a_{i+1} & \alpha_{i+1} \\
\alpha_i & \alpha_{i+1} & \lambda_i \end{pmatrix}$$

If ${a_{i-1}\over a_i}<{a_i\over a_{i+1}}$, then the two
first columns of $M_i$ are linearly independent, and
(\ref{Hodge1_PnxPn_cero}) proves that the last column of $M_i$
is a linear combination of the first two, as wanted.

If instead ${a_{i-1}\over a_i}={a_i\over a_{i+1}}$, then by
(\ref{primera-sustitucion-alphas})
$$\begin{vmatrix} a_i & \alpha_i \\ a_{i+1} &
\alpha_{i+1} \end{vmatrix}=0$$ 
which, together with (\ref{Hodge2_PnxPn_cero}) and
$a_{i+1}\ne0$, implies that the last column of $M_i$ is a
multiple of the middle one.  
\end{pf}

\section{Further examples and remarks}

We start with several examples that illustrate until which
point the method explained in subsection \ref{metodo_general} can
work. In two of them, we will not have any result stating that
subvarieties in the chosen ambient space $Y$ are simply
connected, so that step (ii) of the method cannot be
applied. Anyway, the key point in the method is step (i), in
order to prove that the divisors on subvarieties of $X$ come
numerically from divisors on $Y$. We will not do all the
details, which can be found in the forthcoming PhD thesis
of the second author. We start with a negative example.

\begin{ex}\label{cuadrica}{\rm
Let $X\subset Q_{2n-2}$ be a smooth $n$-dimensional
subvariety of the smooth $(2n-2)$-dimensional quadric. It is
well known that the class of $X$ in the Chow ring of
$Q_{2n-2}$ takes the form $dH^{n-2}$, where $H$ is the hyperplane class.
On the other hand, if $D$ is a (smooth) divisor on $X$, its
class as a cycle in $Q_{2n-2}$ will be
$[D]=\alpha_1A_1+\alpha_2A_2$, where $A_1,A_2$ are the
classes of the two families of $(n-1)$-dimensional linear
spaces contained in $Q_{2n-2}$. Having in mind the
intersection of two linear spaces in $Q_{2n-2}$ depending on
the parity of $n$, equality (\ref{chern_general}) becomes

$$P:=\left\{\begin{matrix}
2\alpha_1\alpha_2-dD^2H^{n-2}=0& $ if
$n$ is even$\\ \\
\alpha_1^2+\alpha_2^2-dD^2H^{n-2}=0& $ if $n$
is odd$\end{matrix}\right.$$
If $S$ is the intersection of $X$ with $n-2$ general
hyperplanes we can thus write (observe that
$D_SH_S=DH^{n-1}=\alpha_1+\alpha_2$)
$$P=\left\{\begin{matrix}
{1\over2}\big((D_SH_S)^2-D_S^2H_S^2\big)
-{1\over2}(\alpha_1-\alpha_2)^2& $if $n$ is even$\\ \\
{1\over2}\big((D_SH_S)^2-D_S^2H_S^2\big)
+{1\over2}(\alpha_1-\alpha_2)^2& $if $n$ is odd.$
\end{matrix}\right.$$
This proves (using the Hodge index theorem for $S$) that, when
$n$ is odd, any divisor on $X$ is numerically a multiple of
the hyperplane section; moreover, since $X$ is a smooth subvariety of
$\P^{2n-1}$, it follows from \cite{BL} that $X$ is simply
connected, and hence we conclude that its Picard group is
isomorphic to $\Z$. However, when $n$ is even, we get that
$D$ is a multiple of $H$ only when $\alpha_1=\alpha_2$, i.e.
$DA_1=DA_2$. Since the equality
$DA_1=DA_2$ defines one linear condition on $Pic X$ (which
could be identically zero), it follows that
$Pic X$ has rank at most two (we need to use again that $X$
is simply connected). The rank-two case can actually occur:
consider the Segre embedding
$\P^1\times\P^{n-1}\hookrightarrow\P^{2n-1}$, whose image can
be defined as the subvariety $X$ consisting of those point
for which the matrix
$$\begin{pmatrix}
X_0&X_1&\ldots&X_{n-1}\\
X_n&X_{n+1}&\ldots&X_{2n-1}
\end{pmatrix}$$
has rank one. Since $n$ is even, it follows that $X$ is
contained in the smooth quadric $Q_{2n-2}$ of equation
$$\begin{vmatrix}X_0 &X_1\\ X_n&X_{n+1}\end{vmatrix}
+\begin{vmatrix}X_2 &X_3\\ X_{n+2}&X_{n+3}\end{vmatrix}
+\ldots+
\begin{vmatrix}X_{n-2} &X_{n-1}\\
X_{2n-2}&X_{2n-1}\end{vmatrix} =0.$$
Strengthening Conjecture 4.5 of \cite{A}, we
suspect that the above example is the only smooth
$n$-dimensional subvariety $X\subset Q_{2n-2}$ whose Picard
group is not generated by $H$.   
}\end{ex}

The next two examples are chosen to see that the method can
keep working even for nonhomogeneous ambient spaces.

\begin{ex}{\rm
Let $Y$ be the blowing up of $\P^6$ at a point. Then the
Picard group (and in fact the whole Chow ring) of $Y$ is
generated by $H$ (the pullback of the hyperplane section of
$\P^6$) and $E$ (the exceptional divisor). Let $X\subset Y$
be a smooth codimension-two subvariety and let $D$ be a
smooth divisor on $X$. We write their classes in the
Chow ring of $Y$ as: 
$$[X]=a_1H^2+a_2E^2$$
$$[D]=\alpha_1H^3+\alpha_2E^3$$

From the intersection table in $X$:
\begin{equation}\label{matriz_expl}
\begin{tabular}{c|cccc}
  & $H^{3}$ & $E^3$ & $DH^2$ & $DE^2$ \\
\hline 
$\begin{matrix} H \\ E \\ D  \end{matrix}$ & 
$\begin{matrix} a_1 \\ 0 \\ \alpha_1  \end{matrix}$ & 
$\begin{matrix} 0 \\ -a_2 \\ -\alpha_2  \end{matrix}$ & 
$\begin{matrix} \alpha_1 \\ 0 \\ \lambda_1  \end{matrix}$ &
$\begin{matrix} 0 \\ -\alpha_2 \\ \lambda_2  \end{matrix}$  
\end{tabular}
\end{equation}
we get that (\ref{chern_general}) takes the form:
$$P:=\alpha_1^2-\alpha_2^2-a_1\lambda_1-a_2\lambda_2=0$$
Then we see that we can get a decomposition
$$P=-{a_1-a_2\over a_1}\begin{vmatrix}
a_1 & \alpha_1  \\ \alpha_1 & \lambda_1
\end{vmatrix}+ {1\over a_1}\begin{vmatrix}
a_1 & 0 & \alpha_1 \\ 0 & -a_2 & -\alpha_2 \\
\alpha_1 & -\alpha_2 & \lambda_1+\lambda_2
\end{vmatrix}$$
The first determinant in this decomposition corresponds
to the intersection matrix of the divisors $H$ and $D$
in the surface obtained as the intersection of $X$ with
two general divisors of class $H$ (this is clearly a
smooth irreducible surface and $H$ is very ample on it),
hence by Hodge index theorem this determinant is
nonpositive. Similarly, the second determinant corresponds
to the intersection matrix of the divisors $H$, $E$ and $D$
in the surface obtained as the intersection of $X$ with
two general divisors of class $H-E$. If we assume for
instance that $H-E$ is ample on $X$, this surface is
smooth irreducible, and hence the determinant is
nonnegative by Hodge index theorem (observe that those three
divisors generate the very ample divisor $2H-E$). Since
$a_1-a_2=(H-E)^4$ (hence strictly positive if we
assume $H-E$ is ample) both determinants
vanish and a simple argument shows that the rank of the
intersection matrix (\ref{matriz_expl}) is two. As we have seen
in the previous sections, this implies (under the
assumption that $H-E$ is ample on $X$) that any divisor
on $X$ depends numerically on $H$ and $E$. 
}\end{ex}

The interest of the next example is that shows that the
method works even if the Picard group of the ambient variety
is not a finitely generated group (and in particular it has
torsion). It also shows (if $C=\P^1$) that the result of
section 3 should remain valid for products of projective
spaces of different dimension.

\begin{ex}{\rm
Set $Y=C\times\P^5$, where $C$ is a smooth curve. The
Picard group of $Y$ is generated by $H$ (the pull-back of
the hyperplane section of $\P^5$) and the pull-back of the
Picard group of $C$ by the projection $p_1:Y\to C$. {}From
the numerical point of view, all the fibers of $p_1$ are
equivalent, and we will denote with $F$ the class of a
fiber. Let $X\subset Y$ be a smooth subvariety of codimension
two and let $D\subset X$ be a smooth divisor. We can write
the numerical classes of $X$ and $D$ as
$$[X]_{num}=a_1H^2+a_2HF$$
$$[D]_{num}=\alpha_1H^3+\alpha_2H^2F$$
We consider now the
intersection matrix:
\begin{equation}\label{matriz_CxP5}
\begin{tabular}{c|cccc}
  & $H_{|X}^{3}$ & $H_{|X}^{2}F_{|X}$ & $DH_{|X}^{2}$ &
$DH_{|X}F_{|X}$ \\
\hline 
$\begin{matrix} H_{|X} \\ F_{X} \\ D  \end{matrix}$ & 
$\begin{matrix} a_2 \\ a_1 \\ \alpha_2  \end{matrix}$ & 
$\begin{matrix} a_1 \\ 0 \\ \alpha_1  \end{matrix}$ & 
$\begin{matrix} \alpha_2 \\ \alpha_1 \\ \lambda_2  \end{matrix}$ &
$\begin{matrix} \alpha_1 \\ 0 \\ \lambda_1  \end{matrix}$  
\end{tabular}
\end{equation}
Equality (\ref{chern_general}) takes now the form:
$$P:=2\alpha_1\alpha_2-a_1\lambda_2-a_2\lambda_1,$$
which can be decomposed as
\begin{equation*}\label{formula_CxP5}
P={1\over a_1}\begin{vmatrix}
a_2 & a_1 & \alpha_2\\ a_1 & 0 & \alpha_1\\ \alpha_2 & \alpha_1 &\lambda_2 \\
\end{vmatrix}-\begin{vmatrix} a_1 &\alpha_1\\ \alpha_1 &
\lambda_1\end{vmatrix}
\end{equation*}
The first determinant in the decomposition corresponds to the
intersection matrix of the classes of $H,F,D$ in the surface
obtained by intersecting $X$ with two divisors of the class
$H$, while the second determinant corresponds to the
intersection matrix of the divisors $H,D$ in the surface
obtained as the intersection of $X$ with a general fiber of
$p_1$ and a hyperplane of $\P^5$. Under the appropriate
assumptions for the irreducibility of those surfaces (even if
$C=\P^1$, Debarre's result cannot provide the irreducibility of the intersection with $HF$ as in section \ref{prod-proj}), we
can apply Hodge index theorem to conclude that both
determinants are zero. In the usual way we conclude that
the rank of the matrix in (\ref{matriz_CxP5}) has rank two
and hence $D$ depends numerically on $H_{|X}$ and
$F_{|X}$. 
}\end{ex}

As these examples (and the rest of the paper) illustrates,
our method works only numerically (and in a mysterious way: we guess that we
are missing something deep that explains when the expression
(\ref{chern_general}) can be decomposed as a sum of nonnegative summands).
Only when there are some topological results (like the ones by Debarre) we
are able to say something about the Picard group, but not as much as we would
expect. In fact we suspect that, in the hypothesis of Theorems
\ref{Barth-Larsen_grass} and \ref{Barth-Larsen_PnxPn} one
could conclude that the Picard group of $X$ is generated
respectively by $H$ and $H_1,H_2$ (and more generally, that
the cohomology coincides in the same range as Barth-Larsen
theorems). 

The fact that analogues of Barth-Larsen theorems seem to
work in general suggests that also there should be a kind of
Hartshorne's conjecture when $2N<3n$. However, the
possible existence of low-rank vector bundles in general
ambient spaces makes difficult to dare with a reasonable
conjecture. For instance, $\G(1,n)$ has a universal
rank-two vector bundle, so that any zero locus of a twist of
it would be subcanonical of codimension two, but not a
complete intersection. In this spirit, the first author
suggested in \cite[Conjecture 4.3]{A} that the only smooth
codimension-two subvarieties in $G(1,n)$ should be complete
intersections or zero loci of twists of the rank-two
universal quotient bundle. In \cite{AF} there is some
evidence for this conjecture, proving it when $n=4$ and the
degree of the variety is at most $25$. 

In arbitrary codimension it seems very risky to state a similar
conjecture for a general ambient space $Y$ of dimension $N$. One could
expect that, under certain general conditions, the restriction map
$H^i(Y,\Z)\to H^i(X,\Z)$ is an isomorphism if $i\le 2n-N$.
Hence, algebraic cycles on $X$ of codimension at most
$\frac{2n-N}{2}$ would lift to cycles on $Y$. Since the top
Chern class of the normal bundle $N$ of $X$ in $Y$ always
lifts to $X$ (by the self-intersection formula) this would
mean that all the Chern classes of $N$ would lift to classes
in $Y$ as long as $N-n-1\le\frac{2n-N}{2}$. Thus a natural
question is to ask under which conditions a smooth subvariety
$X\subset Y$ is the zero locus of a vector bundle of rank
$N-n$ on $Y$ if $3N-2\le4n$ (the only non-trivial
counterexample we know is the Segre embedding of
$\P^1\times\P^{n-1}$ in $Q_{2n-2}$ if $n$ is even). Observe that this
inequality is stronger than the one in Hartshorne's conjecture except for
$N=6$, in which our inequality includes the case $n=4$ (for which
Hartshorne-Serre correspondence gives a positive answer). 

\bigskip

\noindent {\bf Acknowledgements:} We thank Jos\'e Carlos
Sierra and Giorgio Ottaviani for many useful discussions. It
was Jos\'e Carlos Sierra who pointed us the paper
\cite{D}, and helped us to pass our results from numerical equivalence to
linear equivalence.


\begin{thebibliography}{ADSE}

\bibitem{A} E. Arrondo,
\emph{Subcanonicity of codimension two subvarieties}, 
Rev. Mat. Compl. {\bf 18} (2005), 69--80.

\bibitem{AF} E, Arrondo, M.L. Fania, \emph{Evidence to
subcanonicity of codimension two subvarieties of $\G(1,4)$},
to appear on Int. J. Math.

\bibitem{B}  W. Barth,
\emph{Transplanting cohomology classes in complex projective
space}, Amer. J. Math. {\bf 92} (1970), 951--967.

\bibitem{BL} W. Barth, M.E. Larsen, \emph{On the
homotopy-groups of complex projective manifolds}, Math
Scand. {\bf 30} (1972), 88-94.

\bibitem{BV} W. Barth and A. Van de Ven, \emph{On the
geometry in codimension $2$ of Grassmann manifolds}, 
Lecture Notes in Math. {\bf 412},  Springer Verlag (1974),
1--35.

\bibitem{D} O. Debarre, \emph{Th\'eor\`emes de connexit\'e
pour les produits d'espaces projectifs et les
grassmanniennes}, Amer. J. Math. {\bf 118} (1996), no. 6,
1347--1367.

\bibitem{H}  R. Hartshorne, \emph{Varieties of small
codimension in projective space},  Bull. Amer. Math. Soc. 
{\bf 80} (1974), 1017--1032.

\bibitem{KL} S.L. Kleiman, D. Laksov, \emph{Schubert
calculus}, Amer. Math. Monthly {\bf 79} (1972), 1061-1082.

\bibitem{L}  M.E. Larsen, \emph{On the topology of
projective manifolds}, Invent. Math.\/  {\bf 19} (1973),
251--260.

\bibitem{S} A.J. Sommese, \emph{Complex subspaces of
homogeneous complex manifolds. II. Homotopy results},
Nagoya Math. J. {\bf 86}, 101--129. 

\end{thebibliography}
\end{document}